\documentclass[a4paper,11pt]{article}
\usepackage{latexsym}
\usepackage{amsthm}
\usepackage{amsmath}
\usepackage{graphicx}
\usepackage{epsfig}
\usepackage{mathrsfs}
\usepackage{caption}
\usepackage{longtable}
\usepackage[numbers,sort&compress]{natbib}

\date{}

\begin{document}
\title{Enumeration of copermanental graphs\thanks{Supported by NSFC (grant no. 11401044).} }

\author{Shunyi Liu\footnote{Corresponding author.}\\
{\small College of Science, Chang'an University}\\
{\small Xi'an, Shaanxi 710064, P.R. China}\\
{\small E-mail: liu@chd.edu.cn}\\
\\
Jinjun Ren\\
{\small Department of Computer Science, Dingxi Teachers College}\\
{\small Dingxi, Gansu 743000, P.R. China}\\
{\small E-mail: wlkx24206221@163.com}}

\maketitle

\begin{abstract}
Let $G$ be a graph and $A$ the adjacency matrix of $G$. The permanental polynomial of $G$ is defined as $\mathrm{per}(xI-A)$. In this paper some of the results from a numerical study of the permanental polynomials of graphs are presented. We determine the permanental polynomials for all graphs on at most 11 vertices, and count the numbers for which there is at least one other graph with the same permanental polynomial. The data give some indication that the fraction of graphs with a copermanental mate tends to zero as the number of vertices tends to infinity, and show that the permanental polynomial does be better than characteristic polynomial when we use them to characterize graphs.
\end{abstract}
\vspace{2mm}



%%%%%%%%%%%%%%%%%%%%%%%%%%%%%%%%%%%%%%%%%%%%%%%%%%%%%%%%%%%%%
\section{Introduction}

The \emph{permanent} of an $n\times n$ matrix $M$ with entries
$m_{ij}$ $(i,j=1, 2, \dots, n)$ is defined by
\begin{equation*}
\mathrm{per}(M)=\sum_{\sigma}\prod_{i=1}^{n}m_{i\sigma(i)},
\end{equation*}
where the sum is taken over all permutations $\sigma$ of $\{1, 2, \dots, n\}$. This scalar function of the matrix $M$ appears repeatedly in the literature of combinatorics and graph theory in connection with certain enumeration and extremal problems. For example, the permanent of a (0,1)-matrix enumerates perfect matchings in bipartite graphs \cite{LoPl}. The permanent is defined similarly to the determinant. However, no efficient algorithm for computing the permanent is known, while the determinant can be calculated using Gaussian elimination. More precisely, Valiant~\cite{Val} has shown that computing the permanent is $\#$P-complete even when restricted to (0,1)-matrices.

Let $G$ be a graph on $n$ vertices and $A(G)$ the adjacency matrix
of $G$. The \emph{characteristic polynomial} of $G$, $\phi(G,x)$, is defined as the determinant of the characteristic matrix of $A(G)$, i.e.,
\begin{equation*}
\phi(G,x)=\mathrm{det}(xI_n-A(G)),
\end{equation*}
where $I_n$ is the identity matrix of size $n$.

The \emph{permanental polynomial} of $G$, $\pi(G,x)$, is defined as the permanent of the characteristic
matrix of $A(G)$, i.e.,
\begin{equation*}
\pi(G,x)=\mathrm{per}(xI_n-A(G)).
\end{equation*}

It seems that the permanental polynomial was first considered by Turner \cite{Tur}. He in fact considered a graph polynomial called generalized characteristic polynomial which generalizes both the permanental and characteristic polynomials. The permanental polynomials of graphs were first systematically studied by Merris et al. \cite{MeRW}, and the
study of analogous objects in chemical literature was started by Kasum et al. \cite{KaTG}. In these two papers, an important Sachs-form formula was obtained independently which relates the coefficients of the permanental polynomial of a graph with structural properties of the given graph. The literature on permanental polynomial is far less than that on characteristic polynomial (see, for example, \cite{BeFS, Bor, BoJo1, BoJo2, Cas1, Cas2, Che, ChLB, GuCa, HuLB, KaTG, LiTB, MeRW, ToLB, YaZh, ZhLi, ZhLL}). This may be due to the difficulty of computing the permanent.

Recall that two graphs $G$ and $H$ are said to be \emph{cospectral} if they have the same spectrum, and a graph $G$ is \emph{determined} by its spectrum (or equivalently, characteristic polynomial) if any graph cospectral with $G$ is isomorphic to $G$. Analogously, Merris et al. \cite{MeRW} defined two graphs $G$ and $H$ to be \emph{copermanental} if they have the same permanental polynomial. A graph $H$, copermanental but non-isomorphic to a graph $G$, is called a \emph{copermanental mate} of $G$. We say that a graph $G$ is \emph{determined} (or \emph{characterized}) by its permanental polynomial if it has no copermanental mates.

For any graph polynomial, it is of interest to determine its ability to characterize graphs~\cite{Noy}. In \cite{MeRW}, Merris et al. formulated that the permanental polynomial seems a little better than the characteristic polynomial
when it comes to distinguishing graphs which are not trees, since they found that the permanental polynomial distinguishes the five pairs of cospectal graphs of \cite{HaKMR}. It is natural to ask whether the permanental polynomial really performs better than the characteristic polynomial when we use them to distinguish graphs.

Motivated by Merris et al.'s formulation, Liu and Zhang~\cite{LiZh1} studied the copermanental characterizations of some specific graphs. They showed that complete graphs, stars, regular complete bipartite graphs and odd cycles are determined by their permanental polynomials, and found that graphs determined by the characteristic polynomial are not necessarily determined by the permanental polynomial. In particular, it was shown that in general, the paths, even cycles and lollipop graphs cannot be determined by the permanental polynomial~\cite{LiZh1,LiZh2}, while they are determined by their characteristic polynomials~\cite{DaHa1,BoJo,HaLZ}. Recently, C\'{a}mara and Haemers~\cite{CaHa} proved that there is just one pair of cospectral graphs when at most five edges are deleted from $K_n$ (the complete graph on $n$ vertices), whereas Zhang et al.~\cite{ZhWL} showed that all graphs obtained from $K_n$ by removing at most five edges are determined by the permanental polynomial.

The spectral characterizations of graphs have been extensively studied, and it was conjectured that almost all graphs are determined by their characteristic polynomials~\cite{DaHa1,DaHa2}. Some computational results on the spectra of graphs have been obtained. Godsil and McKay~\cite{GoMc} enumerated by computer all graphs on at most 9 vertices, computed their spectra and determined the numbers of graphs for which there exists at least one cospectral mate. Haemers and Spence~\cite{HaSp} extended the computer enumeration to graphs on 10 and 11 vertices. Brouwer and Spence~\cite{BrSp} found the characteristic polynomials for all graphs on 12 vertices, and statistics related to the number of cospectral graphs are obtained. The present data give some indication that possibly almost no graph has a cospectral mate.

The main purpose of this paper is to give the preliminary results of a computational study of copermanentality of graphs. We enumerate by computer the permanental polynomials of all graphs on at most 11 vertices, and determine the number of graphs for which there exists at least one copermanental mate. The numerical data give some indication that the fraction of graphs with a copermanental mate tends to zero as the number of vertices tends to infinity, and show that the permanental polynomial does be better than characteristic polynomial when we use them to characterize graphs.

%%%%%%%%%%%%%%%%%%%%%%%%%%%%%%%%%%%%%%%
\section{Computational Results}
To determine the copermanentality of graphs we first of all have to generate the graphs by computer and then determine their permanental polynomials. All graphs on at most 11 vertices are generated by the well-known nauty and Traces package~\cite{McPi}. Then the permanental polynomials of these graphs are computed by a Maple procedure. Finally we count the number of copermanental graphs.

The results are in Table~\ref{tab:Table 1}. This table lists for $n\le 11$ the total number of graphs on $n$ vertices, the total number of distinct permanental polynomials of such graphs, the number of such graphs with a copermanental mate, the fraction of such graphs with a copermanental mate, and the size of the largest family of copermanental graphs.
\begin{table}[htb]
\scriptsize
\captionsetup{singlelinecheck=off,skip=0pt}
\caption{\small Data of permanental polynomials of graphs on $n\le 11$ vertices}
\label{tab:Table 1}
  \begin{tabular}{rrrrrr}
  \hline
    $n$ & $\#$graphs & $\#$perm. pols & $\#$ with mate & fraction with mate & max. family\\ \hline
     0  & 1          & 1              & 0                      & 0       & 1 \\
     1  & 1          & 1              & 0                      & 0       & 1 \\
     2  & 2          & 2              & 0                      & 0       & 1 \\
     3  & 4          & 4              & 0                      & 0       & 1 \\
     4  & 11         & 11             & 0                      & 0       & 1 \\
     5  & 34         & 34             & 0                      & 0       & 1 \\
     6  & 156        & 153            & 6                      & 0.03846 & 2 \\
     7  & 1044       & 1035           & 17                     & 0.01628 & 3 \\
     8  & 12346      & 12247          & 188                    & 0.01523 & 4 \\
     9  & 274668     & 274153         & 980                    & 0.00357 & 5 \\
     10 & 12005168   & 11999059       & 11869                  & 0.00099 & 7 \\
     11 & 1018997864 & 1018915198     & 163534                 & 0.00016 & 12\\
  \hline
  \end{tabular}
\end{table}

For comparing numerical data directly between the permanental and characteristic polynomials, Table~\ref{tab:Table 2} lists the corresponding data of the characteristic polynomials of graphs on at most 12 vertices~\cite{BrSp,HaSp}.
\begin{table}[htb]
\scriptsize
\captionsetup{singlelinecheck=off,skip=0pt}
\caption{\small Data of characteristic polynomials of graphs on $n\le 12$ vertices}
\label{tab:Table 2}
  \begin{tabular}{rrrrrr}\hline
    $n$ & $\#$graphs    & $\#$char. pols & $\#$ with mate  & fraction with mate & max. family\\ \hline
     0  & 1             & 1              & 0                      & 0 & 1 \\
     1  & 1             & 1              & 0                      & 0 & 1 \\
     2  & 2             & 2              & 0                      & 0 & 1 \\
     3  & 4             & 4              & 0                      & 0 & 1 \\
     4  & 11            & 11             & 0                      & 0 & 1 \\
     5  & 34            & 33             & 2                      & 0.059 & 2 \\
     6  & 156           & 151            & 10                     & 0.064 & 2 \\
     7  & 1044          & 988            & 110                    & 0.105 & 3 \\
     8  & 12346         & 11453          & 1722                   & 0.139 & 4 \\
     9  & 274668        & 247357         & 51039                  & 0.186 & 10 \\
     10 & 12005168      & 10608128       & 2560606                & 0.213 & 21 \\
     11 & 1018997864    & 901029366      & 215331676              & 0.211 & 46 \\
     12 & 165091172592  & 148187993520   & 31067572481            & 0.188 & 128 \\
     \hline
    \end{tabular}
\end{table}

In Table~\ref{tab:Table 1} we see that the fraction of graphs with a copermanental mate appears to be decreasing with the number of vertices. In particular, the fraction decreases sharply at $n=10$, and about 0.099$\%$ of all graphs on 10 vertices are not determined by their permanental polynomials. The data give some indication that possibly the fraction of graphs with a copermanental mate tends to zero as $n$ tends to infinity. An interesting observation from Tables~\ref{tab:Table 1} and \ref{tab:Table 2} is that the fractions of graphs with a copermanental mate are much less than that with a cospectral mate. In addition, the fractions of graphs with a copermanental mate decrease rapidly with the number of vertices, while the fractions of graphs with a cospectral mate increase at first and start to decrease at $n=10$. Clearly, the data show that the permanental polynomial performs better than the charactersitic polynomial when we use them to characterize graphs. However, it is worth pointing out that the study of copermanental characterizations of graphs is more difficult than the cospectral characterizations of graphs (see, for example, \cite{LiZh1,LiZh2,ZhWL}).

Since two graphs with distinct number of edges must have distinct permanental polynomials~\cite{LiZh1}, the enumeration has been carried out for each possible number of edges. The data, differentiated according to the number $m$ of edges, is presented in Tables~\ref{tab:Table_3}--\ref{tab:Table_10}.

\begin{figure}[tphb]
    \begin{center}
\includegraphics[scale=0.6]{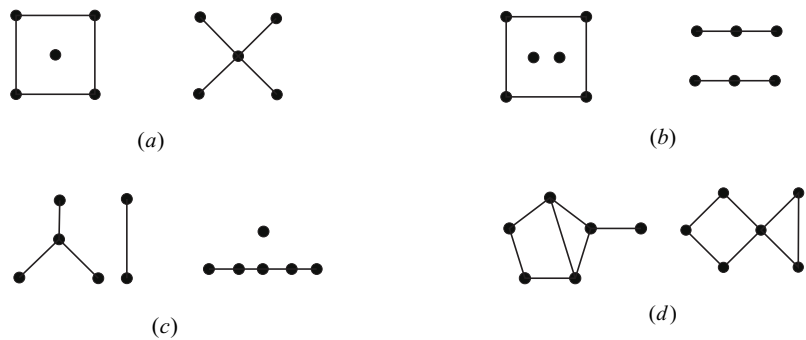}
\end{center}
   \vspace{-5mm}
   \caption{The smallest pairs of cospectral and copermanental graphs.}
   \protect\label{fig:fig_1}
\end{figure}

We end with the smallest pairs of cospectral and copermanental graphs. The smallest pair of cospectral graphs is illustrated in Fig.~\ref{fig:fig_1}(a), and the smallest three pairs of copermanental graphs are illustrated in Fig.~\ref{fig:fig_1}(b), (c) and (d).
\vspace{4mm}

%%%%%%%%%%%%%%%%%%%%%%%%%%%%%%%%%%%%%%%%%%%%%%%%%%%%%%%%%%%%%
\section*{Appendix}

We list the numbers of graphs for all numbers $m$ of edges up to 11 vertices, the numbers of distinct permanental polynomials of such graphs, the numbers of such graphs with a copermanental mate, and the maximum size of a family of copermanental graphs (see Tables~\ref{tab:Table_3}--\ref{tab:Table_10}).
\setlength\LTleft{0pt}
\setlength\LTright{0pt}
\begin{longtable}{@{\extracolsep{\fill}}|r|r|r|r|r|}
\captionsetup{singlelinecheck=off,skip=0pt}
\caption{Graphs on 4 vertices}
\label{tab:Table_3}\\
\hline
$m$ & $\#$graphs & $\#$perm. pols & $\#$ with mate & max. family \\
\hline
\endfirsthead
\multicolumn{5}{l}%
{\tablename\ \thetable\  \textit{(Continued)}} \\
\hline
$m$ & $\#$graphs & $\#$perm. pols & $\#$ with mate & max. family \\
\hline
\endhead
\multicolumn{5}{r}{\textit{(Continued on next page)}} \\
\endfoot
\hline
\endlastfoot
0  & 1          & 1              & 0              & 1 \\
1  & 1          & 1              & 0              & 1 \\
2  & 2          & 2              & 0              & 1 \\
3  & 3          & 3              & 0              & 1 \\
4  & 2          & 2              & 0              & 1 \\
5  & 1          & 1              & 0              & 1 \\
6  & 1          & 1              & 0              & 1 \\
\end{longtable}

\begin{longtable}{@{\extracolsep{\fill}}|r|r|r|r|r|}
\captionsetup{singlelinecheck=off,skip=0pt}
\caption{Graphs on 5 vertices}
\label{tab:Table_4}\\
\hline
$m$ & $\#$graphs & $\#$perm. pols & $\#$ with mate & max. family \\
\hline
\endfirsthead
\multicolumn{5}{l}%
{\tablename\ \thetable\  \textit{(Continued)}} \\
\hline
$m$ & $\#$graphs & $\#$perm. pols & $\#$ with mate & max. family \\
\hline
\endhead
\multicolumn{5}{r}{\textit{(Continued on next page)}} \\
\endfoot
\hline
\endlastfoot
0   & 1          & 1              & 0              & 1 \\
1   & 1          & 1              & 0              & 1 \\
2   & 2          & 2              & 0              & 1 \\
3   & 4          & 4              & 0              & 1 \\
4   & 6          & 6              & 0              & 1 \\
5   & 6          & 6              & 0              & 1 \\
6   & 6          & 6              & 0              & 1 \\
7   & 4          & 4              & 0              & 1 \\
8   & 2          & 2              & 0              & 1 \\
9   & 1          & 1              & 0              & 1 \\
10  & 1          & 1              & 0              & 1 \\
\end{longtable}

\begin{longtable}{@{\extracolsep{\fill}}|r|r|r|r|r|}
\captionsetup{singlelinecheck=off,skip=0pt}
\caption{Graphs on 6 vertices}
\label{tab:Table_5}\\
\hline
$m$ & $\#$graphs & $\#$perm. pols & $\#$ with mate & max. family \\
\hline
\endfirsthead
\multicolumn{5}{l}%
{\tablename\ \thetable\  \textit{(Continued)}} \\
\hline
$m$ & $\#$graphs & $\#$perm. pols & $\#$ with mate & max. family \\
\hline
\endhead
\multicolumn{5}{r}{\textit{(Continued on next page)}} \\
\endfoot
\hline
\endlastfoot
0	& 1	         & 1	          & 0	           & 1 \\
1	& 1	         & 1	          & 0	           & 1 \\
2	& 2	         & 2	          & 0	           & 1 \\
3	& 5	         & 5	          & 0	           & 1 \\
4	& 9	         & 7	          & 4	           & 2 \\
5	& 15	     & 15	          & 0	           & 1 \\
6	& 21	     & 21	          & 0	           & 1 \\
7	& 24	     & 23	          & 2	           & 2 \\
8	& 24	     & 24	          & 0	           & 1 \\
9	& 21	     & 21	          & 0	           & 1 \\
10	& 15	     & 15	          & 0	           & 1 \\
11	& 9	         & 9	          & 0	           & 1 \\
12	& 5	         & 5	          & 0	           & 1 \\
13	& 2	         & 2	          & 0	           & 1 \\
14	& 1	         & 1	          & 0	           & 1 \\
15	& 1	         & 1	          & 0	           & 1 \\
\end{longtable}

\begin{longtable}{@{\extracolsep{\fill}}|r|r|r|r|r|}
\captionsetup{singlelinecheck=off,skip=0pt}
\caption{Graphs on 7 vertices}
\label{tab:Table_6}\\
\hline
$m$ & $\#$graphs & $\#$perm. pols & $\#$ with mate & max. family \\
\hline
\endfirsthead
\multicolumn{5}{l}%
{\tablename\ \thetable\  \textit{(Continued)}} \\
\hline
$m$ & $\#$graphs & $\#$perm. pols & $\#$ with mate & max. family \\
\hline
\endhead
\multicolumn{5}{r}{\textit{(Continued on next page)}} \\
\endfoot
\hline
\endlastfoot
0   & 1	        & 1	             & 0	          & 1 \\
1   & 1	        & 1	             & 0	          & 1 \\
2   & 2	        & 2	             & 0	          & 1 \\
3   & 5	        & 5	             & 0	          & 1 \\
4   & 10	    & 8	             & 4	          & 2 \\
5   & 21	    & 19	         & 4	          & 2 \\
6   & 41	    & 38	         & 6	          & 2 \\
7   & 65	    & 63	         & 3	          & 3 \\
8   & 97	    & 97	         & 0	          & 1 \\
9   & 131	    & 131	         & 0	          & 1 \\
10	& 148	    & 148	         & 0	          & 1 \\
11	& 148	    & 148	         & 0	          & 1 \\
12	& 131	    & 131	         & 0	          & 1 \\
13	& 97	    & 97	         & 0	          & 1 \\
14	& 65	    & 65	         & 0	          & 1 \\
15	& 41	    & 41	         & 0	          & 1 \\
16	& 21	    & 21	         & 0	          & 1 \\
17	& 10	    & 10	         & 0	          & 1 \\
18	& 5	        & 5              & 0	          & 1 \\
19	& 2	        & 2	             & 0	          & 1 \\
20	& 1	        & 1	             & 0	          & 1 \\
21	& 1	        & 1	             & 0	          & 1 \\
\end{longtable}

\begin{longtable}{@{\extracolsep{\fill}}|r|r|r|r|r|}
\captionsetup{singlelinecheck=off,skip=0pt}
\caption{Graphs on 8 vertices}
\label{tab:Table_7}\\
\hline
$m$ & $\#$graphs & $\#$perm. pols & $\#$ with mate & max. family \\
\hline
\endfirsthead
\multicolumn{5}{l}%
{\tablename\ \thetable\  \textit{(Continued)}} \\
\hline
$m$ & $\#$graphs & $\#$perm. pols & $\#$ with mate & max. family \\
\hline
\endhead
\multicolumn{5}{r}{\textit{(Continued on next page)}} \\
\endfoot
\hline
\endlastfoot
0	& 1	        & 1	             & 0	          & 1 \\
1	& 1	        & 1	             & 0	          & 1 \\
2	& 2	        & 2	             & 0	          & 1 \\
3	& 5	        & 5	             & 0	          & 1 \\
4	& 11	    & 9	             & 4	          & 2 \\
5	& 24	    & 20	         & 8	          & 2 \\
6	& 56	    & 48	         & 13	          & 3 \\
7	& 115	    & 102	         & 23	          & 4 \\
8	& 221	    & 200	         & 39	          & 3 \\
9	& 402	    & 392	         & 20	          & 2 \\
10	& 663	    & 652	         & 22	          & 2 \\
11	& 980	    & 971	         & 18	          & 2 \\
12	& 1312	    & 1301	         & 21	          & 3 \\
13	& 1557	    & 1552	         & 10	          & 2 \\
14	& 1646	    & 1643	         & 6	          & 2 \\
15	& 1557	    & 1557	         & 0	          & 1 \\
16	& 1312	    & 1311	         & 2	          & 2 \\
17	& 980	    & 979	         & 2	          & 2 \\
18	& 663	    & 663	         & 0	          & 1 \\
19	& 402	    & 402	         & 0	          & 1 \\
20	& 221	    & 221	         & 0	          & 1 \\
21	& 115	    & 115	         & 0	          & 1 \\
22	& 56	    & 56	         & 0	          & 1 \\
23	& 24	    & 24	         & 0	          & 1 \\
24	& 11	    & 11	         & 0	          & 1 \\
25	& 5	        & 5	             & 0	          & 1 \\
26	& 2	        & 2	             & 0	          & 1 \\
27	& 1	        & 1	             & 0	          & 1 \\
28	& 1	        & 1	             & 0	          & 1 \\
\end{longtable}

\begin{longtable}{@{\extracolsep{\fill}}|r|r|r|r|r|}
\captionsetup{singlelinecheck=off,skip=0pt}
\caption{Graphs on 9 vertices}
\label{tab:Table_8}\\
\hline
$m$ & $\#$graphs & $\#$perm. pols & $\#$ with mate & max. family \\
\hline
\endfirsthead
\multicolumn{5}{l}%
{\tablename\ \thetable\  \textit{(Continued)}} \\
\hline
$m$ & $\#$graphs & $\#$perm. pols & $\#$ with mate & max. family \\
\hline
\endhead
\multicolumn{5}{r}{\textit{(Continued on next page)}} \\
\endfoot
\hline
\endlastfoot
0	& 1	        & 1   	         & 0	          & 1 \\
1	& 1	        & 1	             & 0	          & 1 \\
2	& 2	        & 2	             & 0	          & 1 \\
3	& 5	        & 5	             & 0	          & 1 \\
4	& 11	    & 9	             & 4	          & 2 \\
5	& 25	    & 21	         & 8	          & 2 \\
6	& 63	    & 52	         & 17	          & 4 \\
7	& 148	    & 120	         & 48	          & 5 \\
8	& 345	    & 293	         & 88	          & 5 \\
9	& 771	    & 715	         & 102	          & 3 \\
10	& 1637	    & 1570	         & 127	          & 3 \\
11	& 3252   	& 3210	         & 84	          & 2 \\
12	& 5995	    & 5959	         & 70	          & 3 \\
13	& 10120 	& 10088	         & 64	          & 2 \\
14	& 15615	    & 15574	         & 81	          & 3 \\
15	& 21933	    & 21904	         & 58	          & 2 \\
16	& 27987	    & 27952	         & 69	          & 3 \\
17	& 32403  	& 32376	         & 54	          & 2 \\
18	& 34040	    & 34017	         & 46	          & 2 \\
19	& 32403	    & 32391	         & 24	          & 2 \\
20	& 27987     & 27979	         & 16	          & 2 \\
21	& 21933	    & 21930	         & 6	          & 2 \\
22	& 15615	    & 15610	         & 10	          & 2 \\
23	& 10120	    & 10119	         & 2	          & 2 \\
24	& 5995	    & 5994	         & 2	          & 2 \\
25	& 3252	    & 3252	         & 0	          & 1 \\
26	& 1637	    & 1637	         & 0	          & 1 \\
27	& 771	    & 771	         & 0	          & 1 \\
28	& 345	    & 345	         & 0	          & 1 \\
29	& 148	    & 148	         & 0	          & 1 \\
30	& 63	    & 63	         & 0	          & 1 \\
31	& 25	    & 25	         & 0	          & 1 \\
32	& 11	    & 11	         & 0	          & 1 \\
33	& 5    	    & 5	             & 0	          & 1 \\
34	& 2	        & 2	             & 0	          & 1 \\
35	& 1	        & 1	             & 0	          & 1 \\
36	& 1	        & 1	             & 0	          & 1 \\
\end{longtable}

\begin{longtable}{@{\extracolsep{\fill}}|r|r|r|r|r|}
\captionsetup{singlelinecheck=off,skip=0pt}
\caption{Graphs on 10 vertices}
\label{tab:Table_9}\\
\hline
$m$ & $\#$graphs & $\#$perm. pols & $\#$ with mate & max. family \\
\hline
\endfirsthead
\multicolumn{5}{l}%
{\tablename\ \thetable\  \textit{(Continued)}} \\
\hline
$m$ & $\#$graphs & $\#$perm. pols & $\#$ with mate & max. family \\
\hline
\endhead
\multicolumn{5}{r}{\textit{(Continued on next page)}} \\
\endfoot
\hline
\endlastfoot
0	& 1	        & 1	             & 0	          & 1 \\
1	& 1	        & 1	             & 0	          & 1 \\
2	& 2	        & 2	             & 0	          & 1 \\
3	& 5	        & 5	             & 0	          & 1 \\
4	& 11	    & 9	             & 4	          & 2 \\
5	& 26	    & 22	         & 8	          & 2 \\
6	& 66	    & 53	         & 21	          & 4 \\
7	& 165	    & 129	         & 59	          & 6 \\
8	& 428	    & 330	         & 159	          & 7 \\
9	& 1103	    & 927	         & 306	          & 6 \\
10	& 2769	    & 2434	         & 573	          & 6 \\
11	& 6759	    & 6350	         & 756	          & 6 \\
12	& 15772   	& 15287	         & 921	          & 5 \\
13	& 34663	    & 34288	         & 735	          & 3 \\
14	& 71318	    & 70938	         & 752	          & 3 \\
15	& 136433	& 136083	     & 697	          & 3 \\
16	& 241577	& 241186	     & 777	          & 3 \\
17	& 395166	& 394762	     & 805	          & 3 \\
18	& 596191	& 595739	     & 902	          & 3 \\
19	& 828728	& 828279	     & 896	          & 3 \\
20	& 1061159	& 1060733	     & 852	          & 2 \\
21	& 1251389	& 1251029	     & 718	          & 3 \\
22	& 1358852	& 1358547	     & 610	          & 2 \\
23	& 1358852	& 1358619	     & 466	          & 2 \\
24	& 1251389	& 1251230	     & 318	          & 2 \\
25	& 1061159	& 1061060	     & 198	          & 2 \\
26	& 828728	& 828666	     & 124	          & 2 \\
27	& 596191	& 596150	     & 82	          & 2 \\
28	& 395166	& 395141	     & 50	          & 2 \\
29	& 241577	& 241560	     & 34	          & 2 \\
30	& 136433	& 136421	     & 24	          & 2 \\
31	& 71318	    & 71310	         & 16	          & 2 \\
32	& 34663	    & 34661	         & 4	          & 2 \\
33	& 15772	    & 15771	         & 2	          & 2 \\
34	& 6759	    & 6759	         & 0	          & 1 \\
35	& 2769	    & 2769	         & 0	          & 1 \\
36	& 1103	    & 1103	         & 0	          & 1 \\
37	& 428	    & 428	         & 0	          & 1 \\
38	& 165	    & 165	         & 0	          & 1 \\
39	& 66	    & 66	         & 0	          & 1 \\
40	& 26	    & 26	         & 0	          & 1 \\
41	& 11	    & 11	         & 0	          & 1 \\
42	& 5	        & 5	             & 0	          & 1 \\
43	& 2	        & 2	             & 0	          & 1 \\
44	& 1	        & 1	             & 0	          & 1 \\
45	& 1	        & 1	             & 0	          & 1 \\
\end{longtable}

\setlength\LTleft{0pt}
\setlength\LTright{0pt}
\begin{longtable}{@{\extracolsep{\fill}}|r|r|r|r|r|}
\captionsetup{singlelinecheck=off,skip=0pt}
\caption{Graphs on 11 vertices}
\label{tab:Table_10}\\
\hline
$m$ & $\#$graphs & $\#$perm. pols & $\#$ with mate & max. family \\
\hline
\endfirsthead
\multicolumn{5}{l}%
{\tablename\ \thetable\  \textit{(Continued)}} \\
\hline
$m$ & $\#$graphs & $\#$perm. pols & $\#$ with mate & max. family \\
\hline
\endhead
\multicolumn{5}{r}{\textit{(Continued on next page)}} \\
\endfoot
\hline
\endlastfoot
0   & 1	        & 1	             & 0	          & 1 \\
1	& 1	        & 1	             & 0	          & 1 \\
2	& 2	        & 2	             & 0	          & 1 \\
3	& 5	        & 5	             & 0	          & 1 \\
4	& 11	    & 9	             & 4	          & 2 \\
5	& 26	    & 22	         & 8	          & 2 \\
6	& 67	    & 54	         & 21	          & 4 \\
7	& 172	    & 133	         & 63	          & 6 \\
8	& 467	    & 349	         & 186	          & 8 \\
9	& 1305	    & 1025	         & 453	          & 8 \\
10	& 3664	    & 2981	         & 1110	          & 12 \\
11	& 10250	    & 9017	         & 2110	          & 10 \\
12	& 28259	    & 26102	         & 3867	          & 10 \\
13	& 75415	    & 73121	         & 4348	          & 5 \\
14	& 192788	& 190225	     & 4995	          & 5 \\
15	& 467807	& 465364	     & 4839	          & 4 \\
16	& 1069890	& 1067326	     & 5103	          & 4 \\
17	& 2295898	& 2293123	     & 5539	          & 3 \\
18	& 4609179	& 4605886	     & 6572	          & 4 \\
19	& 8640134	& 8636125	     & 8002	          & 3 \\
20	& 15108047	& 15103213	     & 9655	          & 3 \\
21	& 24630887	& 24625314	     & 11131	      & 4 \\
22	& 37433760	& 37427675	     & 12158	      & 4 \\
23	& 53037356	& 53031054	     & 12596	      & 3 \\
24	& 70065437	& 70059298	     & 12272	      & 3 \\
25	& 86318670	& 86312952	     & 11429	      & 4 \\
26	& 99187806	& 99182803	     & 10004	      & 4 \\
27	& 106321628	& 106317385	     & 8483	          & 3 \\
28	& 106321628	& 106318153	     & 6948	          & 3 \\
29	& 99187806	& 99185067	     & 5476	          & 4 \\
30	& 86318670	& 86316513	     & 4310	          & 4 \\
31	& 70065437	& 70063788	     & 3298	          & 2 \\
32	& 53037356	& 53036062	     & 2588	          & 2 \\
33	& 37433760	& 37432802	     & 1914	          & 4 \\
34	& 24630887	& 24630184	     & 1404	          & 4 \\
35	& 15108047	& 15107546	     & 1002	          & 2 \\
36	& 8640134	& 8639789	     & 690	          & 2 \\
37	& 4609179	& 4608966	     & 426	          & 2 \\
38	& 2295898	& 2295761	     & 274	          & 2 \\
39	& 1069890	& 1069818	     & 144	          & 2 \\
40	& 467807	& 467770	     & 74	          & 2 \\
41	& 192788	& 192774	     & 28	          & 2 \\
42	& 75415	    & 75411	         & 8	          & 2 \\
43	& 28259	    & 28258	         & 2	          & 2 \\
44	& 10250	    & 10250	         & 0	          & 1 \\
45	& 3664	    & 3664	         & 0	          & 1 \\
46	& 1305	    & 1305	         & 0	          & 1 \\
47	& 467	    & 467	         & 0	          & 1 \\
48	& 172	    & 172	         & 0	          & 1 \\
49	& 67	    & 67	         & 0	          & 1 \\
50	& 26	    & 26	         & 0	          & 1 \\
51	& 11	    & 11	         & 0	          & 1 \\
52	& 5	        & 5	             & 0	          & 1 \\
53	& 2	        & 2	             & 0	          & 1 \\
54	& 1	        & 1	             & 0	          & 1 \\
55	& 1	        & 1	             & 0	          & 1 \\
\end{longtable}

%%%%%%%%%%%%%%%%%%%%%%%%%%%%%%%%%%%%%%%%%%%%%%%%%%%%%%%%%%%%%

\end{document}